\RequirePackage{ifpdf}
\ifpdf 
\documentclass[pdftex]{sigma}
\else
\documentclass{sigma}
\fi

\begin{document}

\allowdisplaybreaks

\renewcommand{\thefootnote}{$\star$}

\renewcommand{\PaperNumber}{063}

\FirstPageHeading

\ShortArticleName{Non-Hamiltonian Actions and Lie-Algebra Cohomology of Vector Fields}

\ArticleName{Non-Hamiltonian Actions\\ and Lie-Algebra Cohomology of Vector Fields\footnote{This paper is a
contribution to the Special Issue ``\'Elie Cartan and Dif\/ferential Geometry''. The
full collection is available at
\href{http://www.emis.de/journals/SIGMA/Cartan.html}{http://www.emis.de/journals/SIGMA/Cartan.html}}}

\Author{Roberto FERREIRO P\'EREZ~$^\dag$ and Jaime MU\~{N}OZ MASQU\'E~$^\ddag$}

\AuthorNameForHeading{R. Ferreiro P\'erez and J. Mu\~{n}oz Masqu\'e}

\Address{$^\dag$~Departamento de Econom\'{\i}a Financiera
y Contabilidad I, Facultad de Ciencias Econ\'omicas\\
\hphantom{$^\dag$}~y Empresariales, UCM, Campus de Somosaguas, 28223-Pozuelo de Alarc\'on, Spain}
\EmailD{\href{mailto:roferreiro@ccee.ucm.es}{roferreiro@ccee.ucm.es}}

\Address{$^\ddag$~Instituto de F\'{\i}sica Aplicada, CSIC, C/ Serrano 144, 28006-Madrid, Spain}
\EmailD{\href{mailto:jaime@iec.csic.es}{jaime@iec.csic.es}}

\ArticleDates{Received April 03, 2009, in f\/inal form June 08, 2009;  Published online June 16, 2009}

\Abstract{Two examples of $\mathrm{Dif\/f}^+S^1$-invariant closed
two-forms obtained from forms on jet bundles, which does not admit
equivariant moment maps are presented. The corresponding
cohomological obstruction is computed and shown to coincide with a
nontrivial Lie algebra cohomology class on
$H^2(\mathfrak{X}(S^1))$.}


\Classification{58D15; 17B56; 22E65; 53D20; 53D30; 58A20}

\section{Introduction}

Let $p\colon E\to M$ be a bundle over a compact, oriented
$n$-manifold $M$ without boundary. In~\cite{equiconn}
the forms on the jet bundle of degree greater than the dimension
of the base manifold are interpreted as dif\/ferential forms
on the space of sections by means of the integration map
$\Im \colon \Omega ^{n+k}(J^rE)\to \Omega ^k(\Gamma (E))$
(see Section \ref{inte} for the details).

In particular, if $\alpha$ is a closed $(n+2)$-form on $J^rE$ and
is invariant under the action of a~group~$\mathcal{G}$ on~$E$,
then $\Im \lbrack\alpha ]$ determines a $\mathcal{G}$-invariant
closed two-form on $\Gamma(E)$. In \cite{equiconn} this is applied
to the case of connections on a principal bundle, and in \cite{WP}
to the case of Riemannian metrics. Moreover, in those cases
canonical moment maps for these forms are obtained. The moment
maps are obtained using the fact that the closed two-forms came
from characteristic classes of an invariant connection, and the
moment maps came from the equivariant characteristic classes.

In general, however, if $\alpha \in \Omega ^{n+2}(J^rE)$ is closed
and $\mathcal{G}$-invariant, then $\Im \lbrack \alpha ]$ does not
admit a~moment map necessarily because cohomological obstructions
could exist (see Section \ref{sym} for the details). In this paper
we present two examples where this happens. In the f\/irst example
we consider maps $S^1\to S^1$ and the action of the orientation
preserving dif\/feomorphisms on $S^1$. We def\/ine an invariant closed
two-form on the space $\mathfrak{A} =\left\{ u\colon S^1\to
S^1:\dot{u}(t)\neq 0, \forall t\in S^1 \right\} $, and we show
that the obstruction to the existence of an equivariant moment map
is a nontrivial class in the Lie algebra cohomology of
$\mathfrak{X}(S^1)$.

In the second example we consider regular closed plane curves and
the same group as in the f\/irst example. The invariant cohomology
of the corresponding variational bicomplex is computed in
\cite{AP} and a generator of degree $3$ appears. Again we show
that the closed two-form corresponding to this form does not admit
an equivariant moment map by computing the corresponding
obstruction in the cohomology of the Lie algebra
$\mathfrak{X}(S^1)$ of vector f\/ields on $S^1$.

\section{The integration map}\label{inte}

Let $p\colon E\to M$ be a bundle over a compact,
oriented $n$-manifold $M$ without boundary,
and let $J^rE$ be its $r$-jet bundle
with projections $p_r\colon J^rE\to M$,
$p_{r,s}\colon J^rE\to J^sE$ for $s<r$.
A dif\/feomorphism $\phi \in \mathrm{Dif\/f}E$
is said to be projectable if there exists
$\underline{\phi }\in \mathrm{Dif\/f}M$
satisfying $\phi \circ p=p\circ \underline{\phi}$.
We denote by $\mathrm{Proj}E$ the space
of projectable dif\/feomorphism of $E$,
and we denote by $\mathrm{Proj}^+E$ the subgroup
of elements such that
$\underline{\phi }\in \mathrm{Dif\/f}^+M$, i.e.,
$\underline{\phi }$ is orientation preserving.
The space of projectable vector f\/ields on $E$
is denoted by $\mathrm{proj}E$,
and can be considered as the Lie algebra
of $\mathrm{Proj}E$. We denote by $\phi ^{(r)}$
(resp.\ $X^{(r)}$) the prolongation
of $\phi \in \mathrm{Proj}E$ (resp.\
$X\in \mathrm{proj}E$) to $J^rE$.

Let $\Gamma (E)$ be the space of global sections of $E$ considered
as an inf\/inite dimensional Frechet manifold (e.g. see
\cite[Section I.4]{Hamilton}). For every $s\in\Gamma(E)$ we have
$T_s\Gamma (E)\cong \Gamma (M,s^\ast V(E))$, where $V(E)$ denotes
the vertical bundle of $E$. The group $\mathrm{Proj}E$ acts
naturally on $\Gamma (E)$ in the following way. If $\phi \in
\mathrm{Proj}E$, we def\/ine $\phi _{\Gamma (E)}\in
\mathrm{Dif\/f}\Gamma(E)$ by $\phi _{\Gamma (E)}(s) =\phi \circ
s\circ\underline{\phi }^{-1}$, for all $s\in \Gamma (E)$. In a~similar way, a projectable vector f\/ield $X\in \mathrm{proj}E$
induces a vector f\/ield $X_{\Gamma (E)}\in\mathfrak{X}(\Gamma
(E))$.

Let $\mathrm{j}^r\colon M\times \Gamma(E)\to J^rE$,
$\mathrm{j}^r(x,s)=j_x^rs$ be the evaluation map.
We def\/ine a map
\[
\Im \colon \Omega ^{n+k}(J^rE)\to \Omega ^k(\Gamma (E)),
\]
by setting
\[
\Im [\alpha ]
=\int _M\left( \mathrm{j}^r\right) ^\ast \alpha ,
\]
for $\alpha \in \Omega ^{n+k}(J^rE)$.
If $\alpha \in \Omega ^k(J^rE)$ with $k<n$,
we set $\Im [\alpha ]=0$. The operator $\Im $
satisf\/ies the following properties:

\begin{proposition}[cf.\ \cite{equiconn}]
\label{propF}
For all $\alpha \in \Omega ^{n+k}(J^rE)$
the following formulas hold:
\begin{enumerate}\itemsep=0pt
\item[$1.$]
$\Im \lbrack d\alpha ]=d\Im \lbrack \alpha ]$.

\item[$2.$]
$\Im \lbrack (\phi ^{(r)})^\ast \alpha ]
=\phi _{\Gamma (E)}^\ast\Im \lbrack \alpha ]$,
for every $\phi \in \mathrm{Proj}^+E$.

\item[$3.$]
$\Im \lbrack L_{X^{(r)}}\alpha ]
=L _{X_{\Gamma (E)}}\Im \lbrack \alpha ]$
for every $X\in \mathrm{proj}\,E$.

\item[$4.$]
$\Im \lbrack \iota _{X^{(r)}}\alpha ]
=\iota _{X_{\Gamma (E)}}\Im \lbrack \alpha ]$
for every $X\in \mathrm{proj}\,E$.
\end{enumerate}
\end{proposition}
If $\alpha \in \Omega^{n+k}(J^rE)$,
$s\in \Gamma (E)$, $X_1,\dotsc,X_k
\in T_s\Gamma(E)\cong \Gamma (M,s^\ast V(E))$,
then
\begin{equation}
\Im \lbrack \alpha ]_s(X_1,\dotsc,X_k)
=\int _M(j^rs)^\ast (\iota _{X^{(r)}_k}
\cdots \iota _{X^{(r)}_1}\alpha ).
\end{equation}\label{exp}

There exists a close relationship between
the integration map $\Im $ and the variational
bicomplex, see \cite{localVB} for the details.

More generally, if $R\subset J^rE$ is an open subset and
$\mathcal{R}=\{s\in\Gamma(E):j^r_xs\in R, \forall\, x\in M\}$ is the
space of holomonomic sections of $R$, then the integration map
def\/ines a map $\Im \colon \Omega ^{n+k}(R)\to \Omega
^k(\mathcal{R})$.

\section{Cohomological obstructions
to the existence of moment maps\label{sym}}

Recall the obstructions for the existence of a moment map for an
invariant closed two-form, e.g., see \cite{McS}.

Let $\omega$ be a closed two-form on $M$ and $G$ a group acting on
$M$ and preserving $\omega$. We have an inf\/initesimal action
$\mathfrak{g}\to \mathfrak{X}(M)$, $X\mapsto X_M$.

The action is said to be \emph{weakly Hamiltonian} if for every
$X\in \mathfrak{g}$ the form $\iota _{X_M}\omega $ is exact, i.e.,
if there exists a map $\mu \colon \mathfrak{g}\to C^\infty (M)$
such that for every $X\in \mathfrak{g}$ we have $\iota
_{X_M}\omega =d(\mu (X))$. The action is said to be
\emph{Hamiltonian} if there exists a $G$-equivariant moment map.
At the inf\/initesimal level, if $\mu\colon \mathfrak{g}\to C^\infty
(M)$ is $G$-equivariant then we have $L_X\mu=0$ for every $X\in
\mathfrak{g}$ and the converse is true for connected groups.

If the action is weakly Hamiltonian, the obstructions for the
action to be Hamiltonian lie in $H^2(\mathfrak{g})$: If
$\mu \colon \mathfrak{g}\to C^\infty (M)$ satisf\/ies
$\iota _{X_M}\omega =d(\mu (X))$, we def\/ine
$\tau \colon \mathfrak{g}\times \mathfrak{g}\to \mathbb{R}$ by,
\begin{gather*}
\tau (X,Y) =(L_Y\mu)(X)
 =\mu ([X,Y])+L_{Y_M}(\mu (X)).
\end{gather*}
It can be seen that $\tau$ is closed, that the cohomology class on
$H^2(\mathfrak{g})$ is independent of the $\mu $ chosen, and that
the cohomology class of $\tau$ on $H^2(\mathfrak{g})$ vanishes if
and only if there exists a~moment map $\mu'\colon \mathfrak{g}\to
C^\infty (M)$ such that $\iota _{X_M}\omega =d(\mu' (X))$ and
$L_X\mu'=0$ for every $X\in \mathfrak{g}$. In particular, if the
cohomology class of $\tau$ is not zero, then the action is not
hamiltonian.

\section[Cohomology of smooth vector fields on $S^1$]{Cohomology of smooth vector f\/ields on $\boldsymbol{S^1}$}\label{GyF}

In our examples we apply the preceding results to the action of
the dif\/feomorphism group of $S^1$. Hence the Lie algebra
cohomology of $\mathfrak{X}(S^1)$ appears. This cohomology was
f\/irst computed by Gel'fand and Fuks in \cite{GF} and is well known
(e.g. see \cite{Fuks}). The continuous cohomology
$H(\mathfrak{X}(S^1))$ is isomorphic to the tensor product of a
polynomial ring with one two-dimensional generator $a$ and the
exterior algebra with one three-dimensional generator $b$. The
two-dimensional generator~$a$ is given by,
\begin{gather*}
a(X,Y)  = \int _{S^1}
\left(
\frac{df}{dt}\frac{d^2g}{dt^2}-
\frac{dg}{dt}\frac{d^2f}{dt^2}
\right) dt,
\end{gather*}
where $X=f(t)\frac{d}{dt}$, $Y=g(t)\frac{d}{dt}$.

\section[Mappings $S^1\to S^1$]{Mappings $\boldsymbol{S^1\to S^1}$}

In this f\/irst example we consider the trivial bundle
$E=S^1\times S^1\to S^1$, whose sections are mappings
$u\colon S^1\to S^1$, and the action of $\mathrm{Dif\/f}^+S^1$
on $E$ by $(\phi ,(t,u))\mapsto(\phi (t),u)$
for $\phi\in\mathrm{Dif\/f}^+S^1$ and $(t,u)\in S^1\times S^1$.

The coordinates on $S^1\times S^1$ are denoted by $(t,u)$ and
those on $J^2E$ by $(t,u,\dot{u},\ddot{u})$.

If $X=f(t)\frac{d}{dt}\in\mathfrak{X}(S^1)$,
then its prolongation to $J^2E$ is given by,
\begin{equation}
X^{(2)}=f\frac{\partial }{\partial t}
-\frac{df}{dt}\dot{u}\frac{\partial }{\partial \dot{u}}
-\left(
\frac{d^2f}{dt^2}\dot{u}+2\frac{df}{dt}\ddot{u}
\right)
\frac{\partial}{\partial\ddot{u}}.
\label{X2}
\end{equation}

We consider the open subset $\mathcal{A} \subset J^2E$ def\/ined by
the condition $\dot{u}\neq 0$, which is
$\mathrm{Dif\/f}^+S^1$-invariant, and the form
\[
\sigma =\frac{1}{\dot{u}^2}dt\wedge d\dot{u}
\wedge d\ddot{u}\in\Omega ^3(\mathcal{A}).
\]
It is readily seen that $\sigma $ is closed and that
$L_{X^{(2)}}\sigma =0$ for every $X\in \mathfrak{X}(S^1)$. By
applying the integration operator $\Im $ we obtain a
$\mathrm{Dif\/f}^+S^1$-invariant closed two-form $\omega
=\Im\lbrack\sigma ]\in \Omega ^2(\mathfrak{A})$ on the space
$\mathfrak{A}$ of holonomic sections of $\mathcal{A}$
\[
\mathfrak{A}=\left\{ u\colon S^1\to S^1:\dot{u}(t)\neq 0,\;\forall\,
t\in S^1 \right\} .
\]

In our case, if $s\colon S^1\to S^1\times S^1$
is the section corresponding to $u\colon S^1\to S^1$,
then $T_s\Gamma(\mathcal{A})
\cong \Gamma (S^1,s^{\ast}V(E))
\cong \Gamma(S^1,TS^1)$.

In local coordinates, if
$H=h(t)\frac{\partial }{\partial u}$,
then its prolongation to $J^2E$ is given by,
\[
H^{(2)}=h\frac{\partial }{\partial u}
+\frac{dh}{dt}\frac{\partial }{\partial \dot{u}}
+\frac{d^2h}{dt^2}
\frac{\partial }{\partial \ddot{u}}.
\]
Using this expression and formula \eqref{exp}
we obtain the explicit expression of $\omega $:
\begin{proposition}
If $H,K\in T_s\Gamma(\mathcal{A})$ are given by
$H=h(t)\frac{\partial }{\partial u}$,
$K=k(t)\frac{\partial}{\partial u}$,
and $u\colon S^1 \to S^1$ then
\[
\omega _u(H,K)=\int _{S^1}
\left(
\frac{du}{dt}\right) ^{-2}
\left(
\frac{d^2h}{dt^2}\frac{dk}{dt}
-\frac{dh}{dt}\frac{d^2k}{dt^2}
\right) dt.
\]
\end{proposition}

Moreover, we have $\sigma =d\alpha $, where
$\alpha =\dot{u}^{-1}dt\wedge d\ddot{u}$,
and then,
\begin{gather*}
\iota _{X^{(2)}}\sigma
 =\iota _{X^{(2)}}d\alpha  =L_{X^{(2)}}\alpha
-d(\iota _{X^{(2)}}\alpha),
\end{gather*}
for all $X\in \mathfrak{X}(S^1)$.
By using \eqref{X2} we also obtain,
\[
L_{X^{(2)}}\alpha
=-\frac{d^2f}{dt^2}\frac{1}{\dot{u}}dt\wedge d\dot{u}
=d\left(
-\frac{df}{dt}\frac{d\dot{u}}{\dot{u}}
\right) ,
\]
and hence,
\[
\iota _{X^{(2)}}\sigma
=d\left(
-\frac{df}{dt}\frac{d\dot{u}}{\dot{u}}
-\iota _{X^{(2)}}\alpha
\right)
=d(\rho (X)),
\]
where
\[
\rho (X)
=-\frac{df}{dt}\frac{d\dot{u}}{\dot{u}}
-\iota _{X^{(2)}}\alpha .
\]

Accordingly, the action of $\mathrm{Dif\/f}^+S^1$
on $(\mathfrak{A},\omega )$
is weakly Hamiltonian with moment map
$\mu (X)=\Im \lbrack \rho (X)]$,
$\forall\, X\in \mathfrak{X}(S^1)$,
as we have
\begin{gather*}
\iota _{X_{\mathfrak{A}}}\omega
 =\iota _{X_{\mathfrak{A}}}\Im \lbrack \sigma ]
 =\Im \lbrack \iota _{X^{(2)}}\sigma ]
 =\Im \lbrack d(\rho (X))]
 =d\left(
\Im\lbrack\rho(X)]
\right)
 =d(\mu(X)).
\end{gather*}

The explicit expression of $\mu$ is the following:

\begin{proposition}
If $X=f(t)\frac{d}{dt}$ and $u\colon S^1\to S^1$,
then
\[
\mu (X)_u=-\int _{S^1}
\left(
\frac{du}{dt}
\right)^{-1}
\left(
f\frac{d^3u}{dt^3}
+3\frac{df}{dt}\frac{d^2u}{dt^2}
+\frac{d^2f}{dt^2}\frac{du}{dt}
\right) dt.
\]
\end{proposition}

However, the action is not Hamiltonian.
\begin{proposition}
\label{h3}
If $X=f(t)\frac{d}{dt}$, $Y=g(t)\frac{d}{dt}$, then
\[
\tau(X,Y)=-\int _{S^1}
\left(
\frac{d^2f}{dt^2}\frac{dg}{dt}
-\frac{df}{dt}\frac{d^2g}{dt^2}
\right) dt.
\]
\end{proposition}

\begin{proof}
From the def\/inition of $\tau $ we have
\begin{gather*}
\tau (X,Y)  = L_{Y_{\mathfrak{A}}}\mu (X)
+\mu ([X,Y])
 = \Im \lbrack L_{Y^{(2)}}\rho (X)]
+\Im \lbrack \rho ([X,Y])]\\
\phantom{\tau (X,Y)}{}  = \Im \lbrack L_{Y^{(2)}}\rho X)
+\rho ([X,Y])].
\end{gather*}

From the def\/inition of $\rho $ we have
\begin{gather*}
L_{Y^{(2)}}\rho(X)+\rho ([X,Y])
  =L_{Y^{(2)}}
\left(
-\frac{dx}{dt}\frac{d\dot{u}}{\dot{u}}
\right)
-L_{Y^{(2)}}
\left(
\iota _{X^{(2)}}\alpha
\right)  \\
\phantom{L_{Y^{(2)}}\rho(X)+\rho ([X,Y])  =}{}  -\frac{d}{dt}
\left(
f\frac{dg}{dt}-f\frac{dg}{dt}
\right)
\frac{d\dot{u}}{\dot{u}}
-\iota _{\lbrack X^{(2)},Y^{(2)}]}\alpha .
\end{gather*}
By using \eqref{X2} we obtain
\begin{gather*}
L_{Y^{(2)}}
\left(
\iota _{X^{(2)}}\alpha
\right)
+\iota _{\lbrack X^{(2)},Y^{(2)}]}\alpha
  =\iota _{X^{(2)}}
\left(
L_{Y^{(2)}}\alpha
\right)
 =-f\frac{d^2g}{dt^2}\frac{d\dot{u}}{\dot{u}}
-\frac{df}{dt}\frac{d^2g}{dt^2}dt,\\
L_{Y^{(2)}}
\left(
-\frac{df}{dt}\frac{d\dot{u}}{\dot{u}}
\right)
=-g\frac{d^2f}{dt^2}\frac{d\dot{u}}{\dot{u}}
+\frac{df}{dt}\frac{d^2g}{dt^2}dt,\\
\frac{d}{dt}
\left(
f\frac{dg}{dt}-g\frac{df}{dt}
\right)
\frac{d\dot{u}}{\dot{u}}
=\left(
f\frac{d^2g}{dt^2}-g\frac{d^2f}{dt^2}
\right)
\frac{d\dot{u}}{\dot{u}},
\end{gather*}
and hence,
\[
L_{Y^{(2)}}\rho (X)+\rho([X,Y])
=2\frac{df}{dt}\frac{d^2g}{dt^2}dt.
\]
As
\[
2\frac{df}{dt}\frac{d^2g}{dt^2}dt
=-\left(
\frac{d^2f}{dt^2}\frac
{dg}{dt}-\frac{df}{dt}\frac{d^2g}{dt^2}
\right) dt
+d\left(
\frac{df}{dt}\frac{dg}{dt}
\right) ,
\]
we f\/inally obtain,
\begin{gather*}
\tau (X,Y)
=\int _{S^1}2\frac{df}{dt}\frac{d^2g}{dt^2}dt
=-\int _{S^1}
\left(
\frac{d^2f}{dt^2}\frac{dg}{dt}
-\frac{df}{dt}\frac{d^2g}{dt^2}
\right) dt.\tag*{\qed}
\end{gather*}\renewcommand{\qed}{}
\end{proof}

Accordingly to Section \ref{GyF} the expression
obtained in Proposition \ref{h2} determines
a non-trivial class on the Lie algebra cohomology
of $\mathfrak{X}(S^1)$, and hence $\left\langle
\tau\right\rangle \neq 0$ in $H^2(\mathfrak{X}(S^1))$.
Hence we obtain the following

\begin{corollary}
The action of $\mathrm{Dif\/f}^+S^1$ on $(\mathfrak{A},\omega )$ is
not Hamitonian, i.e., $\omega $ does not admit a~$\mathrm{Dif\/f}^+S^1$-equivariant moment map.
\end{corollary}

\section{Regular plane curves}

Let $E=S^1\times \mathbb{R}^2\to S^1$ be the trivial
bundle.
Global sections of $E$ are none other than mappings
$u\colon S^1\to \mathbb{R}^2$, $u(t)=(x(t),y(t))$.
Coordinates on $E$ are denoted by $(t,x,y)$ and by
$(t,x,y,\dot{x},\dot{y},\ddot{x},\ddot{y})$
the coordinates on $J^2E$.

We consider the open set $\mathcal{R} \subset J^2E$ def\/ined by the
condition $\dot{x}^2+\dot{y}^2\neq 0$, and corresponding to the
$2$-jets of regular curves. The holonomic sections of
$\mathcal{R}$ constitute the space $\mathfrak{Reg}\subset \Gamma
(E)$ of closed regular plane curves,
\[
\mathfrak{Reg}=\left\{
u\colon S^1\to\mathbb{R}^2:
\left\vert
\dot{u}(t)
\right\vert ^2
\neq 0,\; \forall \, t\in S^1
\right\}  .
\]

The group $\mathrm{Dif\/f}S^1$ acts on $E$ by
$(t,u)\mapsto (\phi (t),u)$, for every
$\phi \in \mathrm{Dif\/f}S^1$.
This action induces an action on $J^2E$
and clearly $\mathcal{R}$\ is invariant
under this action. If
$X=f\frac{d}{dt}\in \mathfrak{X}(S^1)$,
then its prolongation to $J^2E$ is given by
\[
X^{(2)}=f\frac{\partial }{\partial t}
-\frac{df}{dt}\dot{x}\frac{\partial }{\partial\dot{x}}
-\frac{df}{dt}\dot{y}\frac{\partial }{\partial\dot{y}}
-\left(
\frac{d^2f}{dt^2}\dot{x}+2\frac{df}{dt}\ddot{x}
\right)
\frac{\partial }{\partial\ddot{x}}
-\left(
\frac{d^2f}{dt^2}\dot{y}+2\frac{df}{dt}\ddot{y}
\right)
\frac{\partial }{\partial \ddot{y}}.
\]
Let us consider the $3$-form on $\mathcal{R}$
given by
\[
\sigma
=v^{-2}dt\wedge dv\wedge d\dot{v},
\]
where $v=\sqrt{\dot{x}^2+\dot{y}^2}$,
$\dot{v}
=v^{-1}(\dot{x}\ddot{x}+\dot{y}\ddot{y})$.
This form appears in \cite{AP} as a generator
of the cohomology for the invariant cohomology
of the variational bicomplex for regular plane
curves. It is easily checked directly that
$L_{X^{(2)}}\sigma=0 $ for every
$X\in\mathfrak{X}(S^1)$ and that $d\sigma =0$.

Let $\omega =\Im \lbrack \sigma ] \in \Omega ^2(\mathfrak{R})$,
where $\Im \colon \Omega ^3(\mathcal{R}) \to \Omega
^2(\mathfrak{Reg})$ is the integration map. By the properties of
$\Im$, we know that $\omega $ is a closed and
$\mathrm{Dif\/f}^+S^1$-invariant two-form on $\mathfrak{Reg}$.

We apply the results of Section \ref{sym} to
the~-- inf\/inite-dimensional~-- case of the
$\mathrm{Dif\/f}^+S^1$-action on $(\mathfrak{Reg},\omega )$.

Let $\alpha \in \Omega ^2(\mathcal{R})$
be the form given by,
$\alpha =v^{-1}dt\wedge d\dot{v}$. Clearly,
we have $d\alpha =\sigma $, and hence
$\iota _{X^{(2)}}\sigma =\iota _{X^{(2)}}d\alpha
=L_{X^{(2)}}\alpha -d(\iota _{X^{(2)}}\alpha )$,
for every $X\in \mathfrak{X}(S^1)$.

As a direct computation shows, we have
\begin{gather}
L_{X^{(2)}}v
=-\frac{df}{dt}v,
\label{f1}\\
L_{X^{(2)}}\dot{v}
 =-\frac{d^2f}{dt^2}v-2\frac{df}{dt}\dot{v},
\label{f2}
\end{gather}
and hence
\begin{gather}\label{Lxa}
L_{X^{(2)}}\alpha
 =-\frac{d^2f}{dt^2}dt\wedge\frac{dv}{v}
 =-d\left(
\frac{df}{dt}
\right)
\wedge \frac{dv}{v}
 =d\left(
-\frac{df}{dt}\frac{dv}{v}
\right) .
\end{gather}
Hence we obtain
\begin{gather*}
\iota _{X^{(2)}}\sigma
  =d\left(
-\frac{df}{dt}\frac{dv}{v}
\right)
-d(\iota _{X^{(2)}}\alpha )
  =d\left(
-\frac{df}{dt}\frac{dv}{v}
-\iota _{X^{(2)}}\alpha
\right) .
\end{gather*}
If we set
\[
\rho (X)=-\frac{df}{dt}\frac{dv}{v}
-\iota _{X^{(2)}}\alpha ,
\]
then the action is weakly Hamiltonian,
with moment map
\[
\mu  (X)=\Im \lbrack \rho  (X)],
\qquad
\forall\, X\in \mathfrak{X}(S^1).
\]

However, the action is not Hamiltonian.

\begin{proposition}
\label{h2}
If $X=f(t)\frac{d}{dt}$,
$Y=g(t)\frac{d}{dt}$,
then
\[
\tau (X,Y)=-\int _{S^1}
\left(
\frac{d^2f}{dt^2}\frac{dg}{dt}
-\frac{df}{dt}\frac{d^2g}{dt^2}
\right) dt.
\]
\end{proposition}

\begin{proof}
We have
\begin{gather*}
\tau (X,Y)
 = L_{Y_{\mathfrak{Reg}}}\mu (X)
+\mu ([X,Y])
 = \Im \lbrack L_{Y^{(2)}}\rho (X)]
+\Im \lbrack \rho ([X,Y])]\\
\phantom{\tau (X,Y)}{}
 = \Im \lbrack L_{Y^{(2)}}\rho (X)
+\rho ([X,Y])],
\end{gather*}
and
\begin{gather*}
L_{Y^{(2)}}\rho (X)
+\rho ([X,Y])
  =L_{Y^{(2)}}
\left(
-\frac{df}{dt}\frac{dv}{v}
\right)
-L_{Y^{(2)}}
\left(
\iota _{X^{(2)}}\alpha
\right) \\
\phantom{L_{Y^{(2)}}\rho (X)
+\rho ([X,Y])  =}{}  -\frac{d}{dt}
\left(
f\frac{dg}{dt}-g\frac{df}{dt}
\right)
\frac{dv}{v}
-\iota _{\lbrack X^{(2)},Y^{(2)}]}\alpha .
\end{gather*}
By using \eqref{f1}, \eqref{f2},
and \eqref{Lxa} we obtain
\begin{gather*}
L_{Y^{(2)}}
\left(
\iota _{X^{(2)}}\alpha
\right)
+\iota _{\lbrack X^{(2)},Y^{(2)}]}\alpha
 =\iota _{X^{(2)}}
\left(
L_{Y^{(2)}}\alpha
\right)  =\frac{d^2g}{dt^2}
\left(
-f\frac{dv}{v}-\frac{df}{dt}dt
\right) ,\\
L_{Y^{(2)}}
\left(
-\frac{df}{dt}\frac{dv}{v}
\right)
  =-\frac{d^2f}{dt^2}g\frac{dv}{v}
+\frac{df}{dt}\frac{d^2g}{dt^2}dt,\\
\frac{d}{dt}
\left(
f\frac{dg}{dt}-g\frac{df}{dt}
\right)
\frac{dv}{v}
 =\left(
f\frac{d^2g}{dt^2}-g\frac{d^2f}{dt^2}
\right)
\frac{dv}{v},
\end{gather*}
and hence,
\[
L_{Y^{(2)}}\rho (X)
+\rho ([X,Y])
=2\frac{df}{dt}\frac{d^2g}{dt^2}dt.
\]
As
\[
2\frac{df}{dt}\frac{d^2g}{dt^2}dt
=-\left(
\frac{d^2f}{dt^2}\frac
{dg}{dt}-\frac{df}{dt}\frac{d^2g}{dt^2}
\right) dt
+d\left(
\frac{df}{dt}\frac{dg}{dt}
\right) ,
\]
we f\/inally obtain
\begin{gather*}
\tau (X,Y)
=\int _{S^1}2\frac{df}{dt}\frac{d^2g}{dt^2}dt
=-\int _{S^1}
\left(
\frac{d^2f}{dt^2}\frac{dg}{dt}
-\frac{df}{dt}\frac{d^2g}{dt^2}
\right) dt.\tag*{\qed}
\end{gather*}\renewcommand{\qed}{}
\end{proof}

As we obtain the same result as that in the preceding example, we
also obtain the following

\begin{corollary}
The action of $\mathrm{Dif\/f}^+S^1$ on $(\mathfrak{Reg},\omega )$
is not Hamitonian, i.e., $\omega $ does not admit a~$\mathrm{Dif\/f}^+S^1$-equivariant moment map.
\end{corollary}

\subsection*{Acknowledgements}

Part of this work started during the stay of the f\/irst author
at Utah State University under the advice of Professor Ian Anderson.
The computations were f\/irst obtained by using \textsc{MAPLE}
package ``Vessiot''. Supported by Ministerio de Ciencia
e Innovaci\'on of Spain under grant \# MTM2008--01386.

\pdfbookmark[1]{References}{ref}
\LastPageEnding

\end{document}